\documentclass[journal]{IEEEtran}

\usepackage{amssymb}
\usepackage[ruled,vlined,linesnumbered]{algorithm2e}
\usepackage{amsthm}
\usepackage{amsfonts}
\usepackage{amsmath}
\usepackage{booktabs}
\usepackage{bm}
\usepackage{color}
\usepackage{cite}
\usepackage{comment}
\usepackage{epsfig}
\usepackage{euscript}
\usepackage{fancybox}
\usepackage{graphicx}
\usepackage{hyperref}
\usepackage{lipsum}
\usepackage{pgf}
\usepackage{psfrag}
\usepackage{pifont}
\usepackage{mathrsfs}
\usepackage{multirow}
\usepackage{setspace}
\usepackage{stfloats}
\usepackage{tabularx}
\usepackage{textcomp}

\newcolumntype{L}[1]{>{\raggedright\arraybackslash}p{#1}}
\newcolumntype{C}[1]{>{\centering\arraybackslash}p{#1}}
\newcolumntype{R}[1]{>{\raggedleft\arraybackslash}p{#1}}

\newcommand{\qa}{{\bf a}}
\newcommand{\qb}{{\bf b}}
\newcommand{\qf}{{\bf f}}
\newcommand{\qs}{{\bf s}}
\newcommand{\qn}{{\bf n}}
\newcommand{\qr}{{\bf r}}
\newcommand{\qv}{{\bf v}}
\newcommand{\qz}{{\bf z}}
\newcommand{\qw}{{\bf w}}
\newcommand{\qu}{{\bf u}}
\newcommand{\qp}{{\bf p}}

\newcommand{\qC}{{\bf C}}
\newcommand{\qL}{{\bf L}}
\newcommand{\qH}{{\bf H}}
\newcommand{\qP}{{\bf P}}
\newcommand{\qQ}{{\bf Q}}
\newcommand{\qR}{{\bf R}}
\newcommand{\qS}{{\bf S}}
\newcommand{\qW}{{\bf W}}
\newcommand{\qlambda}{{\boldsymbol \lambda}}

\newcommand{\bl}[1]{\color{blue}#1}
\newcommand{\rl}[1]{\color{red}#1}

\newcommand*{\argmax}{\operatornamewithlimits{argmax}\limits}

\begin{document}

\title{Local Cyber-Physical Attack for Masking Line Outage and Topology Attack in Smart Grid}

\author{\IEEEauthorblockN{
Hwei-Ming Chung, \IEEEmembership{Student Member, IEEE}, Wen-Tai Li, \IEEEmembership{Member, IEEE}, Chau Yuen, \IEEEmembership{Senior Member, IEEE}, Wei-Ho Chung, \IEEEmembership{Member, IEEE}, Yan Zhang, \IEEEmembership{Senior Member, IEEE}, and Chao-Kai Wen, \IEEEmembership{Senior Member, IEEE} }

\thanks{H.-M. Chung was with Singapore  University  of  Technology  and  Design, Singapore.
He is now with Department of Informatics, University of Oslo, Oslo 0373, Norway, (e-mail: hweiminc@ifi.uio.no).}

\thanks{W.-T. Li and C. Yuen are with the Singapore  University  of  Technology  and  Design,
Singapore 487372 (e-mail: \{wentai\_li, yuenchau\}@sutd.edu.sg).}

\thanks{W.-H. Chung (corresponding author) is with the Department of Electrical Engineering and Institute of Communications Engineering, HsinChu City 300, Taiwan (e-mail: whchung@ee.nthu.edu.tw).}

\thanks{Y. Zhang is with Department of Informatics, University of Oslo, Oslo 0373, Norway (e-mail: yanzhang@ieee.org).}

\thanks{C.-K. Wen is with the Institute of Communications Engineering, National
Sun Yat-sen University, Kaohsiung 804, Taiwan (e-mail: chaokai.wen@mail.nsysu.edu.tw).}

}


\maketitle
\begin{abstract}
Malicious attacks in the power system can eventually result in a large-scale cascade failure if not attended on time.
These attacks, which are traditionally classified into \emph{physical} and \emph{cyber attacks}, can be avoided by using the latest and advanced detection mechanisms.
However, a new threat called \emph{cyber-physical attacks} which jointly target both the physical and cyber layers of the system to interfere the operations of the power grid is more malicious as compared with the traditional attacks. 
In this paper, we propose a new cyber-physical attack strategy where the transmission line is first physically disconnected, and then the line-outage event is masked, such that the control center is misled into detecting as an obvious line outage at a different position in the local area of the power system. 
Therefore, the topology information in the control center is interfered by our attack.
We also propose a novel procedure for selecting vulnerable lines, and analyze the observability of our proposed framework.
Our proposed method can effectively and continuously deceive the control center into detecting fake line-outage positions, and thereby increase the chance of cascade failure because the attention is given to the fake outage. 
The simulation results validate the efficiency of our proposed attack strategy.
\end{abstract}

\section*{Nomenclature} 
\addcontentsline{toc}{section}{Nomenclature}
\begin{IEEEdescription}[\IEEEusemathlabelsep\IEEEsetlabelwidth{$V_1,V_2$}]
\item[A. Sets and Indices]
\item[$j, k$] Bus index.
\item[$l$] Line index.
\item[$n_{b}$]  Total bus number in the system.
\item[$n_{br}$] Total line number in the system.
\item[$i$] Represent $\sqrt{-1}$.
\item[$\mathcal{G}$] Graph representing the topology of the system.
\item[$\mathcal{N}$] Sets of the buses in the system.
\item[$\mathcal{E}$] Sets of the lines in the system.
\item[$\mathcal{L}$] Sets of the buses connected to real line outage position.
\item[$\mathcal{M}$] Sets of the buses connected to fake line outage position.
\item[$\mathcal{B}$] Set of boundary buses in the attack region.
\item[$\mathcal{A}$] Set of buses except boundary buses in the attack region.
\item[$\phi$] Empty set.
\\

\item[B. Variables]
\item[$\qz (\bar{\qz})$]  Measurement vector before (after) attack.
\item[$S_{l}$]  Complex power flow on line $l$. Equal to $P_{l} + i Q_{l}$.
\item[$\overline{S}_{l}$]  Modified complex power flow on line $l$. Equal to $\overline{P}_{l}  + i\overline{Q}_{l} $.
\item[$l_{o}$]  Real outage line.
\item[$l_{m}$]  Fake outage line.
\item[$v_{j}$]  Complex voltage of bus $j$.
\item[$V_{j}$]  Voltage magnitude of bus $j$.
\item[$W_{j}$]  Squared voltage magnitude of bus $j$.
\item[$C_{l}$]  Squared current magnitude of line $l$.
\item[$\theta_{j}$]  Voltage phase of bus $j$.
\item[$P_{j}^{D} (\overline{P}_{j}^{D})$]  Real load of bus $j$ before (after) attack.
\item[$\overline{Q}_{j}^{D}$]  Reactive load of bus $j$ after attack.
\\

\item[C. Parameters]
\item[$P_{e}$]  Parameter error vector of the system.
\item[$\qR$]  Measurement error covariance matrix.
\item[$Y_{l}$]  Admittance of line $l$. Equal to $G_{l} + iB_{l}$.
\item[$Z_{l}$]  Impedance of line $l$.
\item[$\mathbf{L}$]  Line outage distribution factors matrix.
\item[$f_{l}$]  Influence factor of line $l$.
\item[$P_{l}^{max}$]  The thermal limit of line $l$.
\item[$V_{max}$]   Maximum voltage magnitude.
\item[$V_{min}$]   Minimum voltage magnitude.
\\

\item[D. Operators]
\item[$|A|$]   Cardinality of set $A$.
\item[$\qa \cdot \qb$]   Element-wise multiplication of vector $a$ and $b$.
\item[$\Re(a)$]   Real part of complex value $a$.
\\

\item[Other notations are defined in the text.]  

\end{IEEEdescription}

\section{Introduction}
\IEEEPARstart{A} power system plays an important role in supporting modern lives and economy.
Upon a late detection, the initial failures in a power system, may lead to a large-scale cascade failure, and adversely affect the economy and security of a nation\cite{cascade-failure}.
Malicious attacks on a power system can lead to an initially undetectable failure and eventually result in failure if not attended in time.
These attacks can be classified into \emph{cyber} and \emph{physical attacks}.
For \emph{physical attacks}, \cite{2004-salm-terrorist-attack,2007-holm-terrorist-attack} proposed a target selection methods that allow the terrorists to physically attack the power system components (e.g., transmission lines, generators, and transformers), cause a direct power system outage, and triggers cascading failures, e.g. \cite{2014-smith-physical-attack} explored the attacks on the transmission substations in California in $2014$.
However, with the aid of recent technological advancements, power system operators are able to detect these attacks easily and prevent a system failure.

Given that \emph{physical attacks} are easily observable, attackers have resorted to \emph{cyber attacks} where they inject carefully predesigned data to the measurements sent by Supervisory Control and Data Acquisition (SCADA), thereby forcing power system operators into making wrong dispatches.
Although various data processing modules, such as state estimation (SE) and bad data detection, have been built to prevent system operation failures and malicious attacks, these mechanisms can be corrupted by carefully designed cyber attacks.
Accordingly, this topic has attracted much research attention over the past few years \cite{ 2009-liu-FDI,2012-M-FDI,2015-zli-FDI, 2015-Yu-FDI, 2017-FDI-summary, 2015-sankar-mask-outage, 2014-zli-mask-outage,2016-zli-mask-outage-local, 2016_ukra_cyber_phsical, 2016-zhiyi-mask-outage, 2016-zli-mask-outage, 2017-deng-mask-outage-CCPA,2016-sankar-mask-outage,2017-zli-mask-outage-AC,  2013-kim-topology-attack-conf,2013-kim-topology-attack,2016-zli-topology-attack}.

The authors in \cite{2009-liu-FDI} proposed classic false data injection (FDI) attacks that cannot be detected by bad data detection techniques and where the attacker can change the measurements of sensors and capture sufficient information about the power system.
These designed attacks must obey the physical laws, such as Kirchhoff's Current Law, and Kirchhoff's Voltage Law.
The authors in \cite{2012-M-FDI} and \cite{2015-zli-FDI} studied classic FDI attacks with incomplete information about the system, and \cite{2012-M-FDI, 2015-zli-FDI} revealed that cyber attacks can bypass SE and bad data detection techniques even if limited network information is available.
Furthermore, the authors attempted to construct the cyber attack by using principal component analysis without prior information as those in \cite{2015-Yu-FDI}. 
Then, recent development of FDI attacks was summarized in \cite{2017-FDI-summary}.

By targeting both the cyber and physical levels of a power system, the recently emerged \emph{cyber-physical attacks} can interfere with the operations of the system more efficiently than the classic FDI attacks with only pure \emph{cyber attacks}.
\emph{Cyber-physical attacks} can still be classified into line-removing and line-maintaining attacks as described in \cite{2015-sankar-mask-outage}.

In a line-maintaining attack, the attacker physically disconnects the transmission line and simultaneously masks this outage event with altered sensor measurements.
Other forms of advanced line-maintaining attacks have been studied in \cite{2014-zli-mask-outage,2016-zli-mask-outage-local, 2016_ukra_cyber_phsical, 2016-sankar-mask-outage,2016-zhiyi-mask-outage,2017-zli-mask-outage-AC,  2016-zli-mask-outage, 2017-deng-mask-outage-CCPA}.
Specifically, the authors in \cite{2014-zli-mask-outage,2016-zli-mask-outage-local} masked the outage event with a local redistribution attack that was extended to a local attack with incomplete topology information.
A line-maintaining attack was recently launched at the Ukrainian electrical grid in 2015 \cite{2016_ukra_cyber_phsical}, where the physical components of the system were disconnected and the SCADA system was illegally attacked by a third party.
Such attack left $225,000$ customers without power.
The attack was designed based on finding out the line that can cause the most damaging to the system in \cite{2016-zhiyi-mask-outage}.
Then, the authors in \cite{ 2016-zli-mask-outage,2017-deng-mask-outage-CCPA} attempted to modify the PMU data to mask the outage event.
The attack model was further derived using the power flow method \cite{2016-sankar-mask-outage} and SE \cite{2017-zli-mask-outage-AC}.

Meanwhile in a line-removing attack, also called topology attacks, the attacker generates a fake outage event to disturb the regular system operation.
This attack must avoid the trivial solution in order not to be detected by the control center.
Using this approach, the attacker can efficiently mislead the control center with an incorrect network topology and then lead the system to an unstable situation by sending the wrong dispatches.
Line-removing attacks have been studied with both partial and whole information of the system and have been mitigated using countermeasure \cite{2013-kim-topology-attack-conf,2013-kim-topology-attack}.
The authors in \cite{2016-zli-topology-attack} focused on a line-removing attack in a local area and proposed a method for locating the attack region.

Based on the discussions above, the previous approaches have obtained the promising results and demonstrated the potential of the \emph{cyber-physical attacks}.
However, most of these approaches were based on DC model, which is different from the real-world system; also, false data constructed by DC model may cause large residual in AC state estimation \cite{2013-rah-DC-AV-diff}.
Only few studies, such as \cite{2016-sankar-mask-outage,2017-zli-mask-outage-AC,2012-hug-FDI-AC}, proposed the construction of the attack in AC system.
However, in \cite{2016-sankar-mask-outage, 2017-zli-mask-outage-AC}, they still constructed the attack with DC model first and then transformed to AC system.
Then, in \cite{2012-hug-FDI-AC}, it focused more on vulnerability assessment of AC state estimation under cyber-attack.
Moreover, there existed no study that combined line-maintaining and line-removing attacks to create more malicious attack to the power system.
Also, when implementing the \emph{cyber-physical attacks}, one must notice that not all transmission lines in the power system can be selected as attack targets because some lines are strictly protected by the control center.
Only few studies, such as \cite{2015-sankar-mask-outage}, considered the rules of selection.

Motivated by the above observations, we develop a novel attack strategy that combines the line-removing and line-maintaining attack strategies in the AC system.
The attack is implemented in the local area and cannot be easily detected because our design ensures that the physical laws of the power system are satisfied.
Unlike previous studies, we propose a rule for selecting target lines instead of randomly selecting such lines.
We also use the traditional SE method that combines normalized Lagrange multipliers and measurement residuals \cite{2006-Abur-parameter-SE,2016-Abur-parameter-SE,abur-book} to test the effectiveness of our proposed attack strategy.
The contributions of this study are summarized as follows:
\begin{itemize}
\item We propose a novel attack strategy with AC model, which aims to physically attack the transmission line and simultaneously mask the real outage event with the cyber attack by misleading the control center into checking another fake outage line.
With this approach, the control center cannot obtain the accurate information about the topology.
Therefore, the control center needs to develop another method of detecting the system topology.

\item Instead of randomly selecting target lines, we design a target line selection rule based on line outage distribution factor (LODF).
According to the simulation, the proposed target line selection results in higher success rate and no false alarm, which is better than random target line selection.
This method can also help control center to identify the locations vulnerable to attacks as in \cite{2011-sou-secure-index,2016-wen-attack}.

\item We apply the conventional SE and bad data detection techniques to test the effectiveness of the proposed attack. 
The simulation results reveal that the control center detects the fake outage position and that the real outage event is successfully hidden.
This highlights the need of developing another effective detection mechanism.

\end{itemize}

\section{System Model}\label{sec:system_model_problem_formulation}
\begin{figure}[!htb]
\begin{center}
\resizebox{3.5in}{!}{%
\includegraphics*{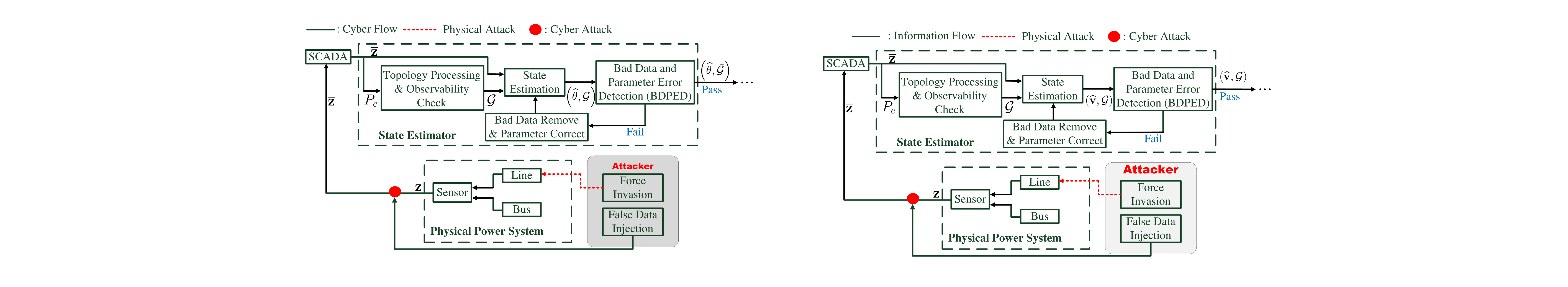} }%
\caption{The system block diagram.}\label{fig:system_model}
\end{center}
\end{figure}

As shown in Fig. \ref{fig:system_model}, the system considered in this study is divided into two parts, namely the SE and the cyber-physical attack model.
We briefly illustrate the SE based on the AC power flow model and the basic calculation of the power system in this section.
Then, our proposed attack strategy is introduced in the next section.
The attack strategy must follow the physical laws introduced in this section.

\subsection{Power Flow Model}\label{subsec:power_flow}
We consider a power transmission network with $n_{b}$ buses and $n_{br}$ lines, and let $\mathcal{N}$ and $\mathcal{E}$ denote the sets of buses and lines, respectively.
The power network can then be represented as a graph denoted as $\mathcal{G}=\{ \mathcal{N}, \mathcal{E}\}$.
Assuming a line $l \in \mathcal{E}$ that connects buses $j$ and $k$, the apparent power of the line flowing from buses $j$ to $k$ denoted as $S_{l}$ can be represented as follows:
\begin{equation}\label{eq:flow_eq}
S_{l} = Y_{l}^{*} |v_{j}|^{2} - Y_{l}^{*}v_{j}v_{k}^{*},
\end{equation}
where $Y_{l}$ is the admittance of line $l$, and $v_{j}$ and $v_{k}$ are the voltages of buses $j$ and $k$ respectively.
The voltage of bus $j$ can be represented in polar form as follows:
\begin{equation}
v_{j} = V_{j} e^{i \theta_{j}},
\end{equation}
where $V_{j}$ is the voltage magnitude, $\theta_{j}$ is the corresponding phase, and $i$ represents $\sqrt{-1}$.
The vector of all power flows is $\qs = [S_{1} \cdots S_{n_{br}}] \in \mathbb{C}^{n_{br} \times 1}$, and the voltage of the buses is denoted as $\qv = [v_{1} \cdots v_{n_{b}}] \in \mathbb{C}^{n_{b}\times 1}$.
The apparent power of line $l$, can be divided into real and reactive power, $P_{l}$ and $Q_{l}$, which are
\begin{align}
& \!\!\! P_{l} = |V_{j}|^{2}G_{l} \!- V_{j} V_{k} \left(  G_{l} \cos (\theta_{j} \!- \theta_{k})  \!+  B_{l} \sin (\theta_{j} \!- \theta_{k}) \right), \label{eq:real_power_eq}\\
& \!\!\! Q_{l} =\!\!  -|V_{j}|^{2}B_{l} \!-\!  V_{j} V_{k} \left(  G_{l} \sin (\theta_{j} \!-\! \theta_{k})  \!-\!  B_{l} \cos (\theta_{j} \!- \theta_{k}) \right), \label{eq:reac_power_eq}
\end{align}
where $G_{l} + iB_{l} = Y_{l}$, and $S_{l} = P_{l} + iQ_{l}$.

\subsection{State Estimation}\label{sec:SE}
Based on the power flow model, we then introduce the state estimator as shown in Fig. \ref{fig:system_model}.
The system states are represented by the voltage of the bus, $\qv$.
Therefore, the measurements received by the SCADA system without attacks can be expressed as follows:
\begin{equation}
\qz = h(P_{e}, \mathcal{G}, \qv) + \qn,
\end{equation}
where $\qz$ usually comprises the measurements of bus injection power and line power flow, and $h(\cdot)$ is the nonlinear function relating to the measurements that depend on the network topology $\mathcal{G}$, network parameter error vector $P_{e}\in \mathbb{C}^{n_{br} \times 1}$ representing the parameter errors, and system state vector $\qv$.
The measurement errors are denoted as $\qn$.
We further denote the measurements modified by the attacker as $\overline{\qz}$.

After obtaining the measurement expression, we adopt the weighted least-squares SE to estimate the system state $\qv$.
The SE problem aims to minimize the sum of squares of the weighted deviations of the measurements estimated from $\overline{\qz}$.
The SE problem is then solved by the following optimization problem where zero parameter errors are assumed:
\begin{subequations} \label{eq:SE_WLS}
\renewcommand\minalignsep{0em}
\begin{align}
\mathcal{F}1 :~  &\min_{\widehat{\qv}, P_{e}}   {\left( \overline{\qz} - h\left(P_{e}, \mathcal{G},   \widehat{\qv}\right)  \right)}^{T} \qR^{-1} {\left( \overline{\qz} - h\left(P_{e}, \mathcal{G}, \widehat{\qv}\right) \right)}  \\
&~\mbox{s.t.}  \quad P_{e} = 0,
\end{align}
\end{subequations}
where $\widehat{\qv}$ is the estimated system state vector, $P_{e}$ is the parameter error vector, and $\qR $ is the measurement error covariance matrix.

\subsection{Bad Data and Parameter Error Detection}\label{sec:bad_data_detec}
After applying SE, we must bypass bad data and parameter error detection to ensure that the measurements are free from bad data or parameter errors.
Accordingly, we apply the normalized residual and Lagrange multiplier method for the detection.

The measurement residual vector can be expressed as:
\begin{equation}
\qr =  \overline{\qz} - h(P_{e}, \mathcal{G}, \qv).
\end{equation}
If the Lagrange multiplier method is applied in (\ref{eq:SE_WLS}), then $\qlambda$ denotes the Lagrange multiplier related to the parameter error.
Given $\qr$ and $\qlambda$, the normalized residual, $\qr^{N}$, and normalized Lagrange multipliers, $\qlambda^{N}$, can be calculated.
The normalized residuals are linked to the corresponding measurements, and the normalized Lagrange multipliers are related to the parameter.
Further details about this procedure can be found in \cite{2006-Abur-parameter-SE,2016-Abur-parameter-SE,abur-book}.
With $\qr^{N}$ and $\qlambda^{N}$, the errors follow the Gaussian distributions, and we choose the largest value between these two parameters of the corresponding line.
If the chosen value is below the identification threshold, then the measurements are free from bad data or parameter errors.
By contrast, the measurement or parameter corresponding to the chosen largest value will be identified as the error.
The part corresponding to this error will be eliminated, and the detection mechanism is reapplied.
Such procedure is performed recursively until no errors are detected.

\section{Proposed Cyber Physical Attack Model }\label{sec:attack_model}
In this section, the attacker block in Fig. \ref{fig:system_model} is illustrated.
We first explain the capabilities of the attackers and the limitations of selecting the target lines, and then introduce the components for launching the attacks separately, namely, selection of target lines, identification of the cyber attack region, and alteration of measurements.
Here, the target lines include the real and the fake outage lines.

\subsection{Introduction of the Attack}\label{subsec:attack_assump}
We assume that the attackers have the following capabilities:
\begin{enumerate}
\item the attacker has knowledge about the topology $\mathcal{G}$ of the entire system;
\item the attacker has the capability to observe the sub-network of $\mathcal{G}$ and perform the power flow calculation for the sub-network; and
\item the attacker only has the capability to change the states of the measurements in the sub-network rather than whole network.
\end{enumerate}
To launch an attack, the attackers are limited to finite sets of target lines because of the following reasons:
\begin{enumerate}
\item the line that connects to a ``three-winding transformer'', or in between two generators cannot be physically attacked;
\item the real and fake outage events cannot take place next to each other; otherwise, the real outage position can be easily observed if the operator goes to repair the line at fake outage position.
Here, the term ``next to'' means that real and fake outage lines have the same end bus.
The mathematical expression can be expressed as 
\begin{equation}
\mathcal{L} \cap \mathcal{M}  \neq \phi;
\end{equation}
\item the generator output cannot be modified;
\item the line injecting power to no-load bus cannot be selected because the load occurring in no load buses can be immediately detected by control center;
\item the load of the buses in the attack region cannot be modified to be negative.
	  Moreover, the difference of the states and measurements before and after the attack must be controlled within a specified range; and
\item if the system is separated into two parts when a line is being attacked, then this line cannot be selected.
\end{enumerate}

\subsection{Real Outage Position Selection for the Physical Attack}\label{sec:line_selection}

When choosing the real outage line, we intend to know the system operation after the outage of a specific line; however, the attacker cannot run the power flow for the whole system.
Therefore, instead of running optimal power flow problem, we employ the LODFs matrix, $\qL \in \mathbb{R}^{{n_{br}}\times {n_{br}}}$, whose definition and calculation can be found in \cite{2009-LODF-cal}.
Specifically, we use the LODF matrix from the DC model to obtain approximate information of a target line if it has been disconnected, because the characteristic of the transmission system is sometimes close to the assumptions of the DC model.
The entry in the $m$-th row and $n$-th column of $\qL$, $l_{m,n}$, represents the fraction of the power flow of the $n$-th line that will be shifted to the $m$-th line when the $n$-th line faces an outage.
By using the LODF matrix, we can define an influence factor, $\qf$, whose $l$-th element can be represented as follows:
\begin{equation}\label{eq:flu_factor}
f_{l} =  {\left(  \left( L_{\{ :,l\}} \right)^{T}  {\rm sign}\left(\Re(\qs) \right) \right)} P_{l},
\end{equation}
where $L_{\{ :,l\}}$ denotes the $l$-th column taken from $\qL$, and $ {\rm sign} \left(\Re(\qs) \right)$ denotes the sign of real power.
The parameter $f_{l}$ represents that the amount of the power flow increment for the whole system if the $l$-th line is disconnected.
In this case, we can determine the real outage line as follows:
\begin{equation}\label{eq:target_line}
l_o= \argmax_l\left\{f_l |~ l=1, \ldots, n_{br}\right\}.
\end{equation}
The $l_{o}$ is the selected real outage position.
Then, the buses connected by $l_{o}$ are assigned to set $\mathcal{L}$.

\subsection{Fake Outage Position Selection and Cyber Attack Region}\label{sec:region_determine}

After presenting the selection of the real outage position, we then illustrate the method of choosing the fake outage position without considering the selection of real outage position.
More clearly, in the part, we assume all the lines remain closed in the system, such that we can ignore the influence from the real outage event.
The idea behind misleading the control center is to let the control center find a fake outage event in the system instead of a real one, thereby hindering the control center from detecting the real outage event and even prompting it to make a wrong operation or decision.
The system faces more risk as the control center spends more time in locating the real outage line.
Moreover, when dispatching the power flows, the control center will avoid assigning the flow to the fake-outage position.
Then, based on this response, the attacker can attempt to create an initial failure.
In the power system, if the power flows are over the thermal limit, meaning the lines are overloaded, this can cause a failure in the power system.
To this end, after the fake outage position is selected, the control center can be misled, and redispatched flow over the residual lines, where the residual lines may end up reaching their thermal limit, and leading to more outage event.
In this context, the optimization problem of selecting the fake outage line is expressed as follows:
\begin{subequations} \label{eq:mislead_opt}
\begin{align}
\mathcal{F}2 : &\max_{w_{l}, \forall l=1,\cdots, n_{br}} 
               &&~ \sum_{l \in \mathcal{E}  } \frac{ \overline{P}_{l} }{ P_{l}^{\max} }  \label{eq:mislead_obj} \\
&\qquad~\mbox{s.t.}   &&~ w_{l} \in \{ 0,1 \},    \label{eq:select_para} \\
&		       &&~ \sum_{l=1}^{n_{br} } w_{l} =  1,  \label{eq:total_select} \\
&		       &&~ \overline{\qp} = \qp +  \left( \qw^{T} \qp \right) \cdot (\qL\qw). \label{eq:power_after_select}
\end{align}
\end{subequations}
Eq.\,(\ref{eq:mislead_obj}) is the objective function that adds the fraction of the real power after certain line has no flow, $\overline{P}_{l}$, to its thermal limit, $P_{l}^{\max}$, for all lines.
Constraints  (\ref{eq:select_para}) and (\ref{eq:total_select}) are the equations related to the selection vector, ${\qw = [w_{1} \cdots w_{n_{br}}]}\in \mathbb{R}^{ n_{br} \times 1} $.
Eq.\,(\ref{eq:power_after_select}) calculates the real power after certain line has no flow, $\overline{\qp}$, based on the LODF matrix.
Therefore, the fake outage position is determined as $l_{m} = \{~l~|~w_l \neq 0,~\forall l=1,\ldots, n_{br} \}$.
The buses linked by the fake outage line, $l_{m}$, are assigned to set $\mathcal{M}$.
Appendix \ref{subsec:proof_line_selection} shows the proof of the selecting criteria.

After the method of selecting the target line is presented, we discuss the method to determine the cyber attack region, based on the limitation that the attackers can only alter the measurement of some selected sensors but not all sensors in the system.
Therefore, we assume that the attacks only have a limited capability to observe and alter a sub-network of $\mathcal{G}$.
To launch an attack, the attacker maliciously changes the measurements in a sub-network of $\mathcal{G}$ denoted as $\mathcal{\overline{G}} = \left\{ \mathcal{\overline{N}}, \mathcal{\overline{E}} \right\}$.
The buses and lines in the attack region are assigned to set $\mathcal{\overline{N}}$ and $\mathcal{\overline{E}}$, respectively.
We further separate set $\mathcal{\overline{N}}$ into sets $\mathcal{A}$ and $\mathcal{B}$.
The boundary buses in $\mathcal{\overline{G}}$ are assigned to set $\mathcal{B}$, when the others are assigned to set $\mathcal{A}$.

The key idea of finding the attack region is that we have to find a new path to re-dispatch the flow, to supply the load of the buses in $\mathcal{M}$, and to obtain a favorable estimate for the power flow of $l_{o}$ and the states of the buses in $\mathcal{L}$.
The sub-network can be obtained using the breadth-first search (BFS) algorithm, which will be discussed later in this paper.

\subsection{Measurements Modification}\label{subsec:modi_measurement}
For the measurement modification, we formulate an optimization problem and meanwhile the physical laws have to be considered as mentioned in the Section \ref{subsec:power_flow}.
To follow these rules, we add the power flow calculation to the constraints.
However, the voltage square in (\ref{eq:flow_eq}) makes the equation become difficult to solve.
The voltage magnitude constraints from the original power flow model also faces the same issue.
To overcome these issues, we use the \textit{DistFlow} model \cite{DistFlow_AC,Hlow_pF_relax}, a convex relaxation technique from the original power flow, to reformulate the AC power flow equations.
In this model, we include two auxiliary variables, namely $\qW = [W_{1} \cdots W_{n_{b}}] \in \mathbb{R}^{n_{b} \times 1}$ and $\qC = [C_{1} \cdots C_{n_{br}}] \in \mathbb{R}^{n_{br} \times 1}$, which denote the squared magnitude of bus voltages and line currents, respectively.
The equation in (\ref{eq:flow_eq}) can then be rewritten as follows:
\begin{equation} \label{eq:Dist_model_flow}
 \lvert S_{l} \rvert^{2} = W_{j} C_{l}.
\end{equation}
Given that the formulation in (\ref{eq:Dist_model_flow}) is still not convex, we apply the following convex relaxation:
\begin{equation} \label{eq:Dist_model_flow_convex}
 \lvert S_{l} \rvert^{2} \leq W_{j} C_{l},
\end{equation}
where (\ref{eq:Dist_model_flow_convex}) is a widely supported second-order cone constraint.
The relation of $\qW$, $\qC$, and $\qs$ can be denoted as
\begin{equation}
 W_{j} - W_{k} = \left(  Z_{l}^{*} S_{l} + Z_{l} S_{l}^{*}  \right) - \lvert Z_{l} \rvert^{2} C_{l} , \\
\end{equation}
where $Z_{l}$ is the impedance of line $l$ that connects buses $j$ and $k$.
The transmission loss of line $l$ can then be calculated as $Z_{l} C_{l}$.

With the power flow model, we then formulate the objective function as minimizing the required attacker ability, that is, we minimize the difference in the measurement before and after the attack. 
These measurements may include the voltage magnitudes, voltage phases, loads of buses, and power flows of the lines.
However, the power flows of the lines are closely related to the voltage magnitudes, voltage phases, and loads of buses.
Therefore, the objective function can be defined as follows:
\begin{equation} 
J =  \| \overline{\qs} - \qs  \|_{2},
\end{equation}
where $\overline{\qs}$ is the modified power flow in the attack region.

According to the discussions above, the optimization can be formulated with $J$ as the objective function.
Then, power flow equations and altered range of the measurements are regarded as the constraints.
The optimization formulation is then formulated as
\begin{subequations} \label{eq:AC_PF}
\renewcommand\minalignsep{-1.5em}
\begin{align}
\mathcal{F}3 : &\min_{ \qW, \overline{\qP}^{D}, \overline{\qQ}^{D},\overline{\qs} } \qquad J     \label{eq:AC_total_obj}\\
&\mbox{s.t.}~ \notag \\
          & W_{j} = |v_{j}|^{2}, & & \forall j \in \mathcal{B}, \label{eq:AC_voltage_boundary} \\
		   & V_{\min}^{2} \leq  W_{j} \leq V_{\max}^{2} ,  & & \forall j \in \mathcal{A}, \label{eq:AC_voltage_inside}\\
	       &  W_{j} \!\!-\!\! W_{k} \! = (  Z_{l}^{*} \overline{S}_{l} \!+\! Z_{l} \overline{S}_{l}^{*} ) \!\!-\!\! \lvert Z_{l} \rvert^{2} C_{l}, & & \forall j, k \!\in \mathcal{\overline{N}}, l \!\in \mathcal{\overline{E}}, \label{eq:AC_voltage_relation} \\
		   & (1-\tau) P_{j}^{D} < \overline{P}_{j}^{D} < (1+\tau) P_{j}^{D},     & &  \forall j \in \mathcal{\overline{N}},  \label{eq:AC_load_range}\\
	       & \overline{P}_{j}^{D} + i \overline{Q}_{j}^{D}  = \sum_{l}  \left( \overline{S}_{l} + Z_{l} C_{l} \right),  & & \forall j \in A , l \in \mathcal{E}, \label{eq:AC_power_injec}\\
			& -P_{l}^{\max} < \overline{P}_{l}  < P_{l}^{\max} , & &  \forall l \in \mathcal{\overline{E}}, \label{eq:AC_line_limit} \\
			& \left|\overline{S}_{l} \right|^{2}\leq W_{j} C_{l},  &  &  \forall j \in \mathcal{\overline{N}}, l \in \mathcal{\overline{E}}. \label{eq:AC_power_cal}
\end{align}
\end{subequations}
Eq.\,(\ref{eq:AC_voltage_boundary}) shows that the voltage magnitude of the boundary buses must remain the same, and Eq.\,(\ref{eq:AC_voltage_inside}) shows that the voltage magnitudes in $\mathcal{A}$ must be controlled within a specified range.
$V_{\max}$ and $V_{\min}$ represent the upper and lower bound of the voltage magnitude, respectively.
Eq.\,(\ref{eq:AC_voltage_relation}) shows relation of $\qW$, $\qC$, and the modified apparent power at the ``from" end of line $l$,  $\overline{S}_{l}$, that connects buses $j$ and $k$.
Eq.\,(\ref{eq:AC_load_range}) shows that the real load of bus $j$ inside the region before, $P_{j}^{D}$, and after modification, $\overline{P}_{j}^{D}$, must be maintained within a specified range.
Moreover, $\tau$ indicates the modification range.
The power injected into the bus must meet the load listed in Eq.\,(\ref{eq:AC_power_injec}) and $\overline{Q}_{j}^{D}$ is the reactive load of bus $j$.
Eq.\,(\ref{eq:AC_line_limit}) shows that the modified flows of the $l$-th line, $\overline{P}_{l}$, must be managed under the thermal limits.
Eq.\,(\ref{eq:AC_power_cal}) calculates the apparent power flow in the attack region.

If the attacker wishes to implement the attack using DC model, then the optimization problem can be easily formulated.
The DC model assumes that if all voltage magnitudes are set to $1$, then no transmission loss takes place, and the phase difference can be neglected.
Based on this assumption, the formulation is explained in detail in \cite{2017-Chung-combine-attack}.

\section{Implementation Strategy}\label{sec:modi_imple}
\begin{figure}
\begin{center}
\resizebox{3.2in}{!}{%
\includegraphics*{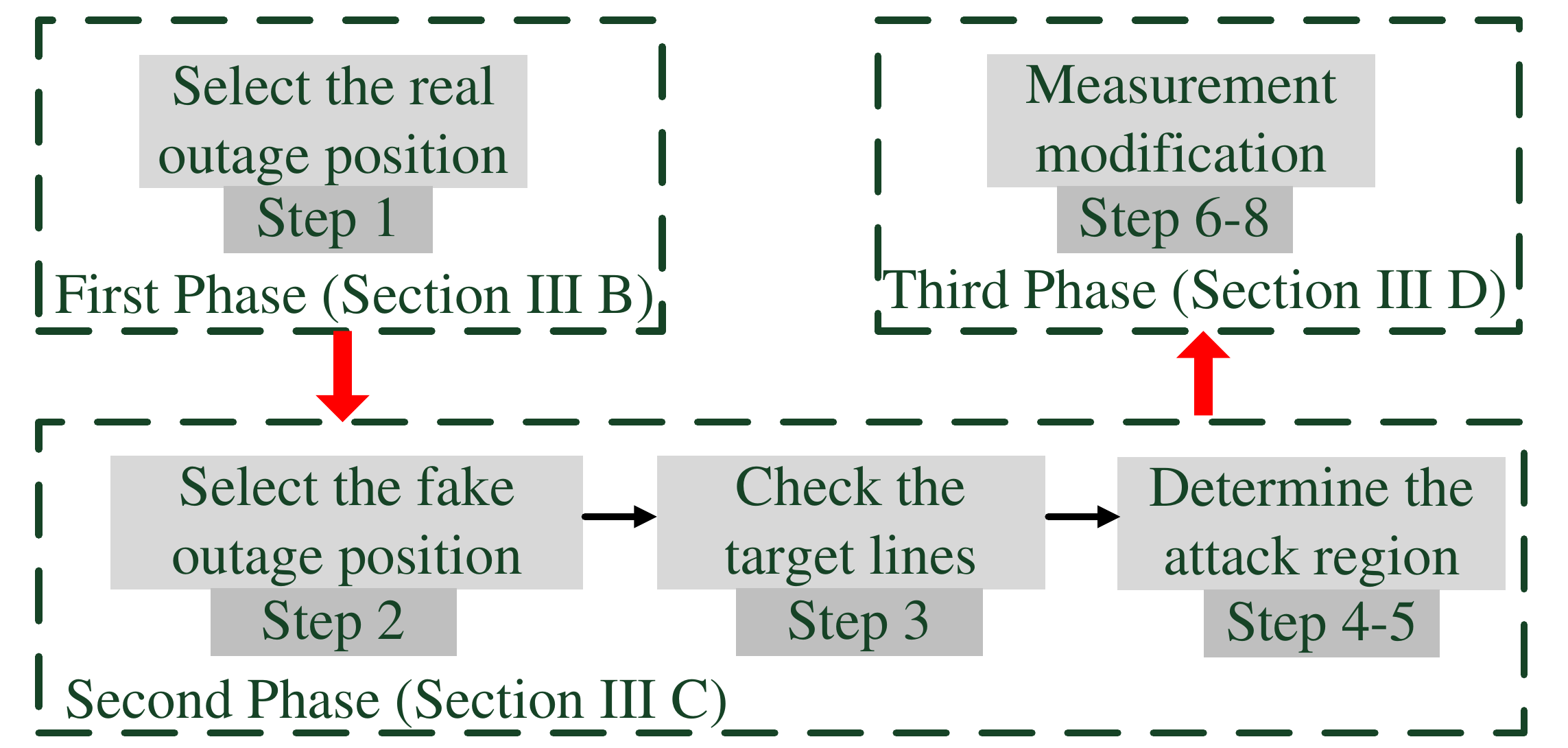} }%
\caption{The implementation strategy of the proposed attack.}\label{fig:attack_procedure}
\end{center}
\end{figure}

In this section, we explain the implementation strategy of the proposed attack scheme by gathering the attack components presented in Section \ref{sec:attack_model}.
This strategy is divided into three phases including total $8$ steps as shown in Fig. \ref{fig:attack_procedure}.
Step $1$ is the first phase which is used to find the real outage position.
Then, the second phase is from Step $2$ to $5$ containing the fake outage position selection and attack region determination.
In the final, the third phase focuses on the measurement modification listed in Step $6$ to $8$.

In the first phase, we use Eq.\,(\ref{eq:target_line}) to determine the real outage position.
The procedure is outlined as follows:

Step 1: We select the line whose disconnection generates the greatest influence on the system.
Therefore, $l_{o}$ is selected as the description in (\ref{eq:target_line}), and the buses linked by $l_{o}$ are assigned to $\mathcal{L}$.
We then use LODF to calculate the power flow after $l_{0}$ faces an outage.

With the real outage position obtained in Step 1, we then use Problem $\mathcal{F}2$ to select the corresponding fake outage position.
After determining the target lines, we must check if the selection fulfills the rules described in Section \ref{subsec:attack_assump}.
We then detailly describe the steps of deciding the attack region.
The region of sub-network $\mathcal{\overline{G}}$ is obtained by using the BFS algorithm to find the shortest path for redispatching the power flow.
The procedure is outlined as follows:

Step 2: We apply the fake outage position selection algorithm listed in Algorithm \ref{ago:line_sel} to select the fake outage position.
We construct vector $\qu$ at the beginning of the algorithm, and then assign each element of $\qu$ with the objective function of Problem $\mathcal{F}2$ on the basis of the results of the exhaustive search.
Line $l_{0}$ must avoid being selected, when line $l_{m}$ obtains the largest value of the objective function.

Step 3: After selecting the target lines, we must check if these lines are reasonable or follow the rules described in Section \ref{subsec:attack_assump}.
If these lines are not reasonable or do not follow such rules, we then eliminate the unreasonable line $l_{o}$ from $\qf$ or $l_{m}$ from $\qu$, and then start again from Step 1.
Otherwise, we proceed to the Step 4.

Step 4: Assume that $l_{m}$ flows from buses $j$ to $k$.
The trivial solution is that we only add and subtract the flow amount to and from buses $j$ and $k$, respectively.
However, this trivial solution can be easily recognized by the control center.
In this case, we only add the flow amount to the load of bus $j$ and try to find another path to supply the load at bus $k$ to prevent the application of a trivial solution.

Step 5: To find a path for supplying bus $k$, we initially consider $\mathcal{\overline{N}}$ and $\mathcal{\overline{E}}$ in $\mathcal{\overline{G}}$ as empty sets.
Afterward, we use BFS algorithm described in Algorithm \ref{ago:BFS_alo} to find the shortest path for redispatching the flow.
The path obtained from Algorithm \ref{ago:BFS_alo} is regarded as the sub-network $\mathcal{\overline{G}}$.
We add $l_{o}$ to $\mathcal{\overline{E}}$ and the buses in $\mathcal{L}$ to $\mathcal{\overline{N}}$ as the attack region.
The sub-network $\mathcal{\overline{G}}$ is therefore determined.

\begin{algorithm}
\caption{Fake outage position selection algorithm}
\label{ago:line_sel}
 \DontPrintSemicolon
\KwIn{Power flow $\qp$, LODF matrix $\qL$}
\KwOut{fake outage position $l_{m}$}
$\qu = [u_{1} \cdots u_{n_{br}}]\in \mathbb{R}^{ n_{br} \times 1}  $ \;
\For{l  = 1 \KwTo $n_{br}$ }{
\eIf{ $l = l_{o}$ }
{
$u_{l}=0$\;
}
{ $\overline{\qp} = \qp +  P_{l} \cdot \qL_{ \{ 1:n_{br}, l  \} }$ \;
 $ u_{l} = \sum_{l \in \mathcal{E} \setminus  l_{o} } \frac{ \overline{P}_{l} }{ P_{l}^{\max} } $
}
$l_{m} = \argmax_{l} \left \{ u_{l} | l=1, \ldots, n_{br}  \right\}$
   }
\end{algorithm}

From step 1 to step 3, the attack strategy involves many condition checks, and therefore the attackers may spend most of time on searching the target lines. 
However, if the attackers can analyze the power system and find out the lines, which do not follow the limitation in Section \ref{subsec:attack_assump}, before launching the attacks, the target lines selection can converge more quickly.
Then, with the target lines and the attack region, the modification is then based on the solution to Problem $\mathcal{F}3$.
The procedure is outlined as follows:

Step 6: After selecting the attack region, we first set the power flow of the fake outage position to $0$.
Then, solve the optimization Problem $\mathcal{F}3$ based on the AC model listed in (\ref{eq:AC_PF}).
However, if the attacker wishes to employ the DC model, then the formulation in \cite{2017-Chung-combine-attack} must be solved.
This formulation is a convex optimization problem that can be dealt with using many existing algorithms and toolboxes.

Step 7: If the problem has no solution, then the current attack region cannot satisfy the constraints.
We then reapply Algorithm \ref{ago:BFS_alo}, and go back to Step 6.
Otherwise, we proceed to Step 8.

Step 8: We set $\overline{\qz} = \qz$ and replace the power flow measurements of $\overline{\qz}$ in $\overline{\mathcal{G}}$ with the $\overline{\qs}$ from the solution to Problem $\mathcal{F}3$.

With the proposed mechanism, the normalized Lagrange multiplier of fake outage position will have the largest value compared to others.
In this case, the operator can easily detect an outage event happening at the fake position.
Appendix \ref{subsec:proof_effect} presents the proof of this statement.

\begin{algorithm}
\caption{BFS algorithm for finding fake outage position}
\label{ago:BFS_alo}
 \DontPrintSemicolon
\KwIn{System topology $\mathcal{G}$, bus $k$, sub-network $\mathcal{\overline{G}}$, line $l_{o}$}
\KwOut{Sub-network $\mathcal{\overline{G}}$}
Find a bus $g$ having a generator and is the nearest to bus $k$\;
Current system configuration is $\mathcal{W} = \left\{ \mathcal{N}, \mathcal{E}\setminus \left\{ \mathcal{\overline{E}}, l_{o} \right\} \right\}$. \;
$g$ : starting bus.   $k$ : destination bus.  \;
let the bus $g$ be the {\it progress bus} and the level $t=0$.
Rest buses are set as {\it unvisited buses}.\;
Search all of the {\it unvisited buses} connected to the buses in {\it progress buses}.
Put such {\it unvisited buses} to {\it progress buses} and previous {\it progress buses} are assigned as {\it visited buses}.
\;
\eIf{ $j \in$ {\it progress bus} }
{
go to step 11 of Algorithm \ref{ago:BFS_alo}. \;
}
{repeat step $5$ of Algorithm \ref{ago:BFS_alo} again.\;
$t = t + 1$. \;
}
Backtrack from the destination bus to the starting bus level-by-level, and identify the shortest path.
The buses and lines in the path are given to $\mathcal{\overline{N}}$ and $\mathcal{\overline{E}}$ respectively.\;
\end{algorithm}

\section{Case Study}\label{sec:simulation}

\begin{figure}
\begin{center}
\resizebox{2.8in}{!}{%
\includegraphics*{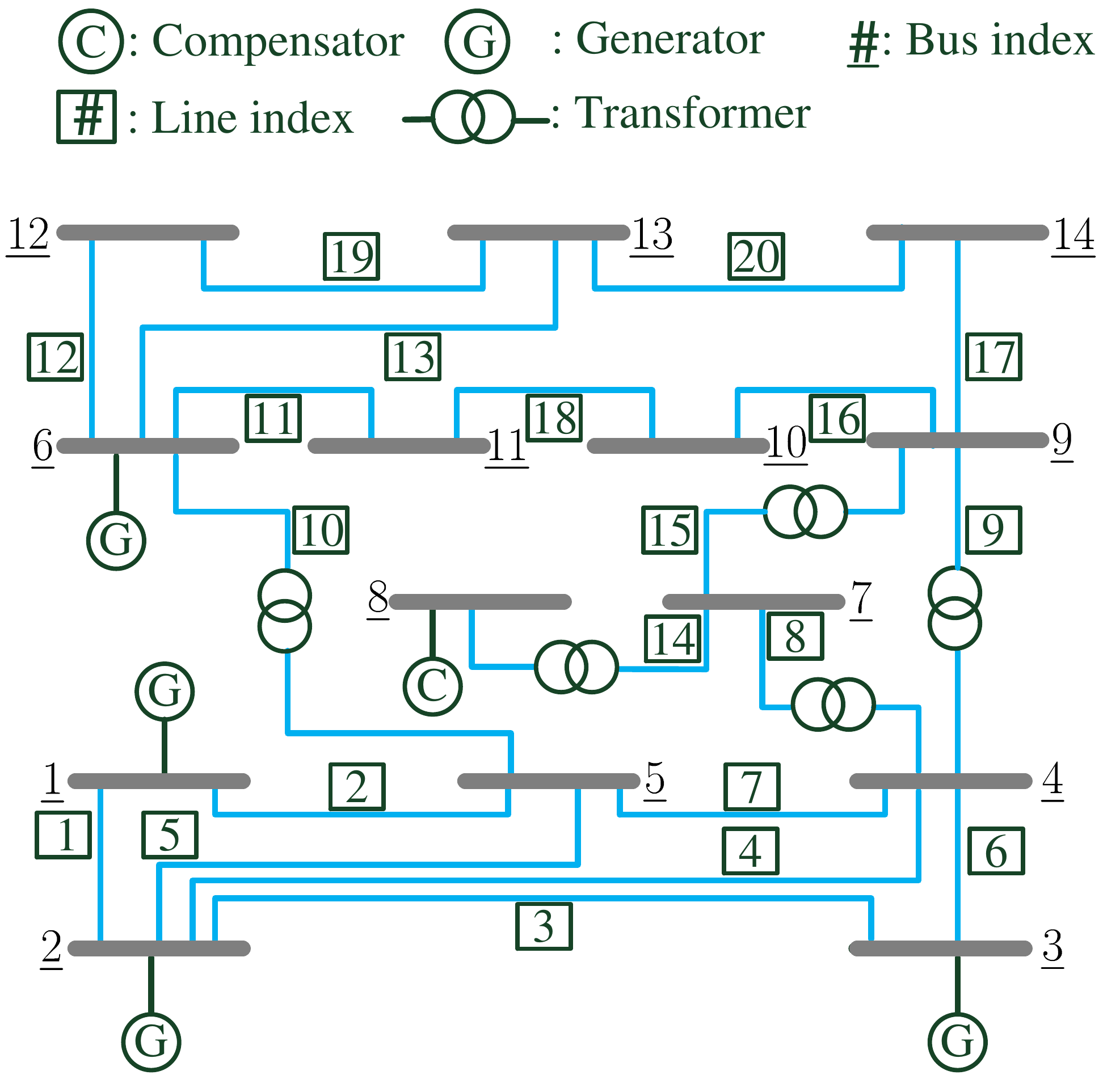} }%
\caption{IEEE 14-bus test system \cite{IEEE_14bus}.} \label{fig:14bus_system}
\end{center}
\end{figure}

In this section, we adopt the IEEE 14-bus \cite{IEEE_14bus}, 24-bus \cite{IEEE_24bus}, and 118-bus systems \cite{IEEE_118bus} to validate our proposed attacking mechanism.
Specifically, we employ the 14-bus test system to illustrate the proposed method in detail.
Fig. \ref{fig:14bus_system} shows the topology of the 14-bus system, when Table \ref{tb:flow_limit} lists the thermal flow limit of each line.
Table \ref{tb:flow_limit} lists the thermal flow limit of each line in 14-bus system.
We then use 118-bus test system to demonstrate the efficiency of the proposed attack mechanism in a large system, and the corresponding thermal limits can be found in \cite{IEEE_118bus}.
Without any specification, the modification rage, $\tau$, for all measurements and loads is set to $25 \%$.
The voltage magnitude must be controlled between $1.05$ p.u. and $0.95$ p.u..
The measurements used for the bad data and parameter error detection include the real power of all lines at the  ``from'' end, the voltage magnitude of all buses, and the voltage phase of the reference bus.
The errors for all measurements are assumed to be $n_{i} \sim N(0, 0.001)$.
The identification threshold of bad data and parameter error detection is set to $3$ which is outside the $99.80 \%$ confidence interval.
The blue and red colors in the simulation results denotes the real and fake outage positions, respectively.
We use the software toolbox MATPOWER \cite{MATPOWER} to run the power flow and provide the initial information of the system.
To solve the optimization problem in (\ref{eq:AC_PF}), we use CVX \cite{CVX}, a package for solving convex programs.
The SE results are obtained through Gauss-Newton method\cite{abur-book}.
Also, all the simulations for computation were conducted with an Intel i7-6700 computer with 3.4 GHz CPU and 8GB RAM.
The increment ratio of the $j$-th measurement, $\Delta_{j}$, which denotes the incremental amount compared to the original measurement, can be expressed as follows:
\begin{equation}
\Delta_{j} = \frac{|\bar{z}_{j}| - |z_{j}|}{ |z_{j}| } \times 100\% .
\end{equation}

\begin{table} \small
\begin{center}
\caption{The thermal flow limit and influence factor of IEEE 14-bus system}\label{tb:flow_limit}
\begin{tabular}{|c|c|c|c|c|c|c|c|c|c|c|c|}
\hline
Line  number  		& 1     	& 2   		& 3			& 4    		& 5  	 \\
\hline
Limit (MW)  		& $200$     & $100$     & $100$      & $100$ 	& $100$  \\
\hline
$f_{l}$    			& $-28.37$     & $14.60$    	& $53.23$      & $4.81$    	& $-32.18$   \\
\hline \hline
Line  number  		& 6     	& 7   		& 8			& 9    		& 10  	 \\
\hline
Limit (MW)  		& $50$     & $100$     & $50$      & $100$ 		& $100$  \\
\hline
$f_{l}$    			& $-10.71$     & $-30.43$    	& $-4.75$      & $6.86$    	& $10.65$   \\
\hline \hline

Line  number  		& 11     	& 12   		& 13		& 14    	& 15  	 \\
\hline
Limit (MW)  		& $50$     & $20$   	& $50$      & $50$	 	& $50$  \\
\hline
$f_{l}$    			& $2.16$     & $-7.33$    	& $12.45$      & NaN    	& -33.05   \\
\hline \hline

Line  number  		& 16     	& 17   		& 18		& 19   		& 20  	 \\
\hline
Limit (MW)  		& $20$     & $20$     & $20$      & $20$ 		& $20$  \\
\hline
$f_{l}$    			& $-2.62$     & $-0.22$    	& $-1.36$      & $-1.36$    	& $1.10$   \\
\hline 
\end{tabular}
\end{center}
\end{table}

\subsection{Implementing the Attack in the 14 Bus System}\label{subsec:14bus_implement}

\begin{figure}
\begin{center}
\resizebox{3in}{!}{%
\includegraphics*{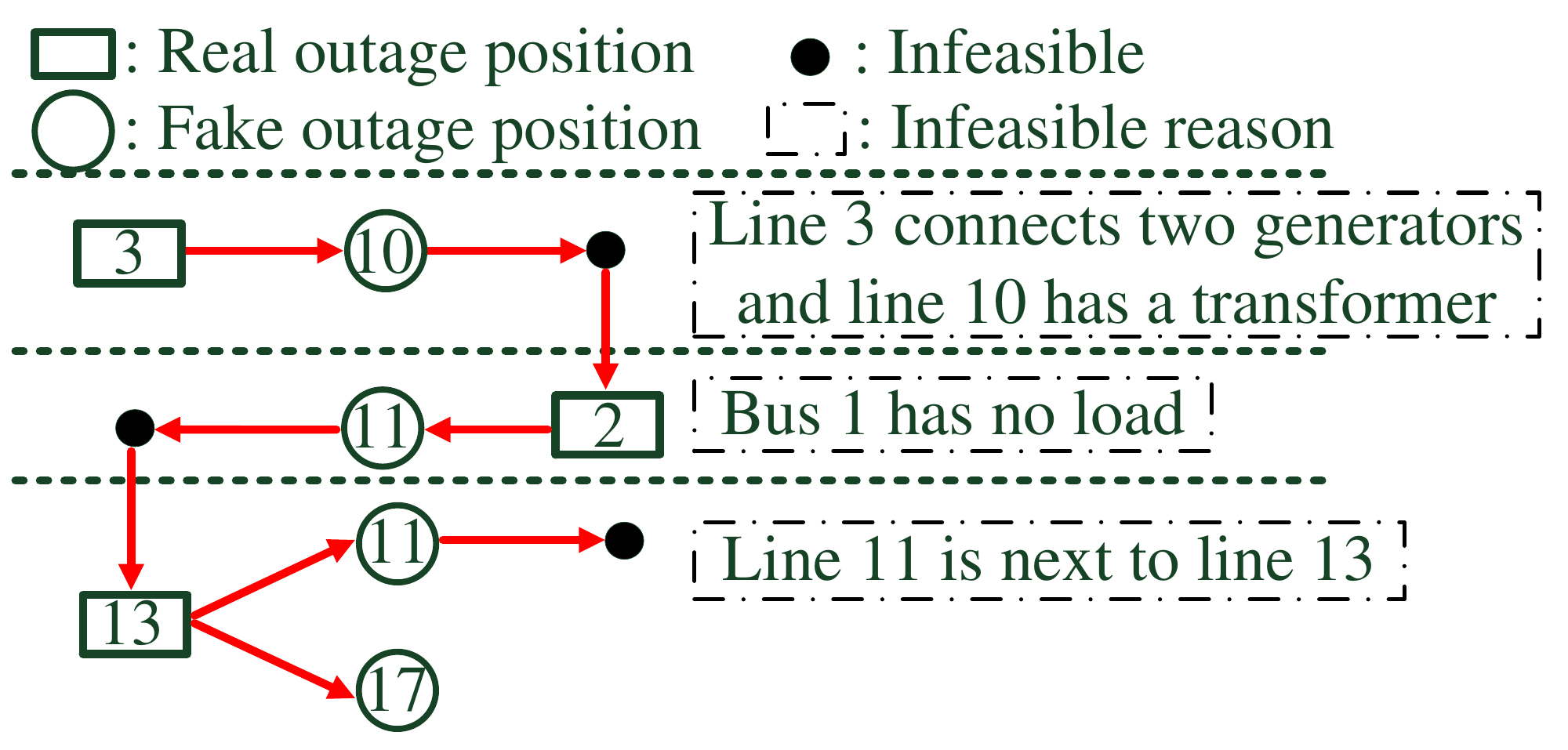} }%
\caption{Target lines selection procedure of 14-bus system.}\label{fig:line_select_14bus}
\end{center}
\end{figure}

We select the target lines following Steps 1 to 3 outlined in Section \ref{sec:modi_imple}.
The influence factors of all the lines are listed in Table \ref{tb:flow_limit}.
Line $3$ has the largest value in the influence factor and the corresponding fake outage position is line $10$ based on the solution of Problem $\mathcal{F}2$.
However, line $3$ is connected to two generators and a transformer is located at line $10$.
Therefore, these target lines are not feasible, and we have to find another sets of real and fake outage positions.
The second largest value in Table \ref{tb:flow_limit} is line $2$. 
Line $2$ is selected as real outage position, and line $11$ is the corresponding fake outage position.
However, line $2$ connects bus $1$ and bus $5$, and bus $1$ has no load; line $2$ cannot be the real outage line.
We then selected the third largest value in Table \ref{tb:flow_limit}.
Line $13$ is selected as the real outage position here, and line $11$ is still the corresponding fake outage position.
However, line $11$ is next to line $13$.
We then eliminate line $11$ from the possible solutions, and resolve the Problem $\mathcal{F}2$ again.
Following the proposed recursive way, we eventually select lines $13$ and $17$ as the real and fake outage lines, respectively.
Fig. \ref{fig:line_select_14bus} illustrates the target line selection process in detail.

Given that the flowing path of the fake outage position is from buses $9$ to $14$, we must find a path for supplying the load of bus $14$.
The nearest generator is located at bus $6$.
Therefore, we use Algorithm \ref{ago:BFS_alo} to find the shortest path from the starting bus, bus $6$, to the destination bus, bus $14$.
Table \ref{tb:set_detail} summarizes the attack region based on the results of Algorithm \ref{ago:BFS_alo}.
Table \ref{tb:modi_load_phase} and \ref{tb:modi_flow} present the measurements before and after modification based on the solutions of Problem $\mathcal{F}3$.
The computation time for solving Problem $\mathcal{F}3$ is $0.40$ second.

\begin{table} \small
\begin{center}
\caption{The sets used in the modification for 14-bus system}\label{tb:set_detail}
\begin{tabular}{|c|c|c|c|c|c|c|c|c|}
\hline
Set   	&  Bus number		& Description  \\
\hline
$\mathcal{A}$    & $12, 13, 14$  & The buses in the attack region          \\
\hline
$\mathcal{B}$    & $6$           & The boundary bus of the attack region \\
\hline
$\mathcal{L}$     & $6, 13$      & The buses connecting the real outage line \\
\hline
$\mathcal{M}$    & $9, 14$       & The buses connecting the fake outage line       \\
\hline
\hline
Set   	&  Line number		&　Description  \\
\hline
$\mathcal{\overline{E}}$    & $12, 13, 19, 20$       & The lines in the attack region      \\
\hline
\end{tabular}
\end{center}
\end{table}

\begin{table*}
\begin{center}
\caption{The voltage and load before and after attack with 14-bus test system}\label{tb:modi_load_phase}
\begin{tabular}{|c|c|c|c|c|c|c|c|c|c|c|}
\hline
Bus number   	& \multicolumn{3}{c|}{Voltage Magnitude (p.u.)}  & \multicolumn{3}{c|}{Voltage Phase (Angle)}		  &  \multicolumn{3}{c|}{Load (MVA)	}  \\
\hline
			    & Before 	& After 	& Increment 			& Before 		  & After		& Increment 		   		& Before 					& After 					& Increment \\
\hline
$6$     		& $1.0500$  & $1.0500$  & $\phantom{0}0.00\%$	& $-0.1499$       & $-0.1499$	& $\phantom{0}0.00\%$     	& $11.20 + 7.50i$ 			& $ 8.96 + 7.65i $    		& $-12.60\%$ \\
\hline
$12$    		& $1.0350$  & $1.0363$  & $\phantom{0}0.13\%$	& $-0.1648$       & $-0.1656$	& $\phantom{0}0.48\%$       & $\phantom{0}6.10 + 1.60i$ & $ 4.58 + 0.95i $    		& $-25.83\%$ \\
\hline
$13$ 	        & $1.0304$  & $1.0290$  & $-0.14\%$	& $-0.1658$       & $-0.1682$	& $\phantom{0}1.48\%$      & $13.50 + 5.80i$        	&  $ \!\!\! 10.13 + 5.66i $ & $-21.02\%$   \\
\hline
$14$    		& $1.0161$  & $1.0042$ & $-1.17\%$ 	& $-0.1778$       & $-0.2033$	& $13.78\%$      			& $14.90 + 5.00i$        	& $\!\!\! 11.18 + 1.49i $  	& $-28.24\%$  \\
\hline

\end{tabular}
\end{center}
\end{table*}

We then perform SE as well as bad data detection by using these modified measurements.
We perform the detection by calculating the normalized residual and Lagrange multipliers, and then sort the results shown in Table \ref{tb:14bus_bad_data_detec}(a) in a descending order.
Table \ref{tb:14bus_bad_data_detec}(a) shows two large Lagrange multipliers that are related to $x_{17}$ and $x_{20}$ and are larger than the identification threshold.
We then eliminate those measurements that are related to $x_{17}$ and $x_{20}$ before reapplying bad data detection.
Table \ref{tb:14bus_bad_data_detec}(b) shows the results of the second round of detection.
The largest value in Table \ref{tb:14bus_bad_data_detec}(b) is much lower than the threshold.
Therefore, we successfully mislead the control center into detecting an error on the fake outage position.

\begin{table} \small
\begin{center}
\caption{The power flow before and after attack for 14-bus system }\label{tb:modi_flow}
\begin{tabular}{|c|c|c|c|c|c|c|c|c|}
\hline
Line number   	& \multicolumn{3}{c|}{Power flow (MVA)	} 		  \\
\hline
	    & Before 						& After						 & Increment \\
\hline
$12$    & $\phantom{0}7.58 + 2.58i$     & $\phantom{0}7.79 + 2.32i$  & $\phantom{-00}1.51\%$    \\
\hline
${\bl 13}$    & $16.90 + 7.36i$      	& $18.93 + 7.47i$      		 & $\phantom{-0}10.40\%$  \\
\hline
${\rl 17}$    & $10.49 + 3.41i$      	& $0$         				 & $-100.00\%$\\
\hline
$19$    & $\phantom{0}1.41 + 0.83i$     & $\phantom{0}2.85 + 0.62i$  & $\phantom{-0}78.26\%$       \\
\hline
$20$    & $\phantom{0}4.60 + 1.98i$     & $11.39 + 1.92i$         	 & $\phantom{-}130.64\%$ \\
\hline

\end{tabular}
\end{center}
\end{table}

\begin{table} 
\begin{center}
\caption{The detection results for 14-bus system}\label{tb:14bus_bad_data_detec}
\begin{tabular}{|c|c|c|c|c|c|c|}
   \multicolumn{2}{c}{(a) First Round}                                                          &  \multicolumn{3}{c}{ ~~~~(b) Second Round}     \\
   \cline{1-2} \cline{4-5}
   \multirow{1}{*}{Parameter}   & \multirow{2}{*}{$\lambda_{j}^{N}$, $r_{j}^{N}$ }  & & \multirow{1}{*}{Parameter}   & \multirow{2}{*}{$\lambda_{j}^{N}$, $r_{j}^{N}$ } \\
   \multirow{1}{*}{Measurement} &                                                   & & \multirow{1}{*}{Measurement} &　\\ \cline{1-2} \cline{4-5}   \cline{1-2} \cline{4-5}
   ${\rl x_{17}}, x_{20}$ 		& $3.4160$                                          & & $x_{12}, {\bl x_{13}}, x_{19}$ 	& $0.8050$ \\    \cline{1-2} \cline{4-5}
   ${\bl x_{13}}$				& $3.1934$                                          & & $x_{7}$ 		          		& $0.0430$ \\    \cline{1-2} \cline{4-5}
   $x_{11}, x_{16}, x_{18}$		& $2.6841$                                          & & $x_{5}$ 						& $0.0340$ \\    \cline{1-2} \cline{4-5}
   $x_{12}, x_{19}$				& $0.4600$                                          & & $r_{12}, r_{13}, r_{19}, r_{33}$& $0.0255$ \\    \cline{1-2} \cline{4-5}
 \cline{1-2} \cline{4-5}
\end{tabular}
\end{center}
\end{table}

\subsection{Consequences of the Attack in the 118 Bus System}\label{subsec:118bus_implement}
Following the same procedure, we apply the proposed attack mechanism to the 118-bus test system.
We then select lines $4$ and $21$ as the positions of the real and fake outage events, respectively.
The attack region contains $16$ buses and $19$ lines.
We also employ the modified measurements for the bad data and parameter error detection, which results are summarized in Table \ref{tb:118bus_bad_data_detec}.
In the simulation, the computation time for solving Problem $\mathcal{F}3$ is $0.80$ second.
For the first round detection presented in Table \ref{tb:118bus_bad_data_detec}(a), the operator easily observes that line $21$ has an obvious error.
After eliminating the measurements related to line $21$, Table \ref{tb:118bus_bad_data_detec}(b) shows no error in the measurements.
Therefore, in the large system, the operator can also be misled into detecting a fake outage event.

\begin{table} 
\begin{center}
\caption{The detection results for 118-bus system}\label{tb:118bus_bad_data_detec}
\begin{tabular}{|c|c|c|c|c|c|c|}
   \multicolumn{2}{c}{(a) First Round}                                                          &  \multicolumn{3}{c}{ ~~~~(b) Second Round}     \\
   \cline{1-2} \cline{4-5}
   \multirow{1}{*}{Parameter}   & \multirow{2}{*}{$\lambda_{j}^{N}$, $r_{j}^{N}$ }  & & \multirow{1}{*}{Parameter}   & \multirow{2}{*}{$\lambda_{j}^{N}$, $r_{j}^{N}$ } \\
   \multirow{1}{*}{Measurement} &                                                   & & \multirow{1}{*}{Measurement} &　\\ \cline{1-2} \cline{4-5}   \cline{1-2} \cline{4-5}
   ${\rl x_{21}}$	 			& $10.9228$                                         & & $x_{8}, x_{37}$	    	    	& $2.2111$ \\    \cline{1-2} \cline{4-5}
   $x_{23}, x_{24}$				& $\phantom{0}7.4409$                               & & $x_{20}, x_{22}$           		& $2.2018$ \\    \cline{1-2} \cline{4-5}
   $x_{26}$						& $\phantom{0}7.0079$                               & & $x_{96}$   						& $1.9532$ \\    \cline{1-2} \cline{4-5}
   $x_{36}$						& $\phantom{0}6.0390$                               & & $x_{30}$           				& $1.8193$ \\    \cline{1-2} \cline{4-5}
 \cline{1-2} \cline{4-5}
\end{tabular}
\end{center}
\end{table}

\subsection{Observability and Effect of Target Lines Selection} \label{subsec:compare_random_selection}
The bad data and parameter error detection can be influenced by noise.
We collect the results of $1,000$ Monte Carlo simulations and study the influence of target line selection.
Specifically, we examine the following cases:
\begin{itemize}
\item {\it Cases $1$, $3$, and $5$ } : Select target lines based on the proposed method in the 14-bus, 24-bus, and 118-bus systems, respectively; and
\item {\it Cases $2$, $4$, and $6$ } : Randomly select target lines in the 14-bus, 24-bus, and 118-bus systems, respectively.
\end{itemize}
For the 24-bus system, $\tau$ is set to $60\%$ because each line has a very huge power flow.
Therefore, the measurements in a local area with a small $\tau$ cannot be easily modified.
If the loads vary over a certain range, it cannot pass the sanity check from the control center even the results cannot be detected by the bad data detection.
This setting just helps us demonstrate the proposed attack strategy in 24-bus system.

\begin{figure}
\begin{center}
\resizebox{3.2in}{!}{%
\includegraphics*{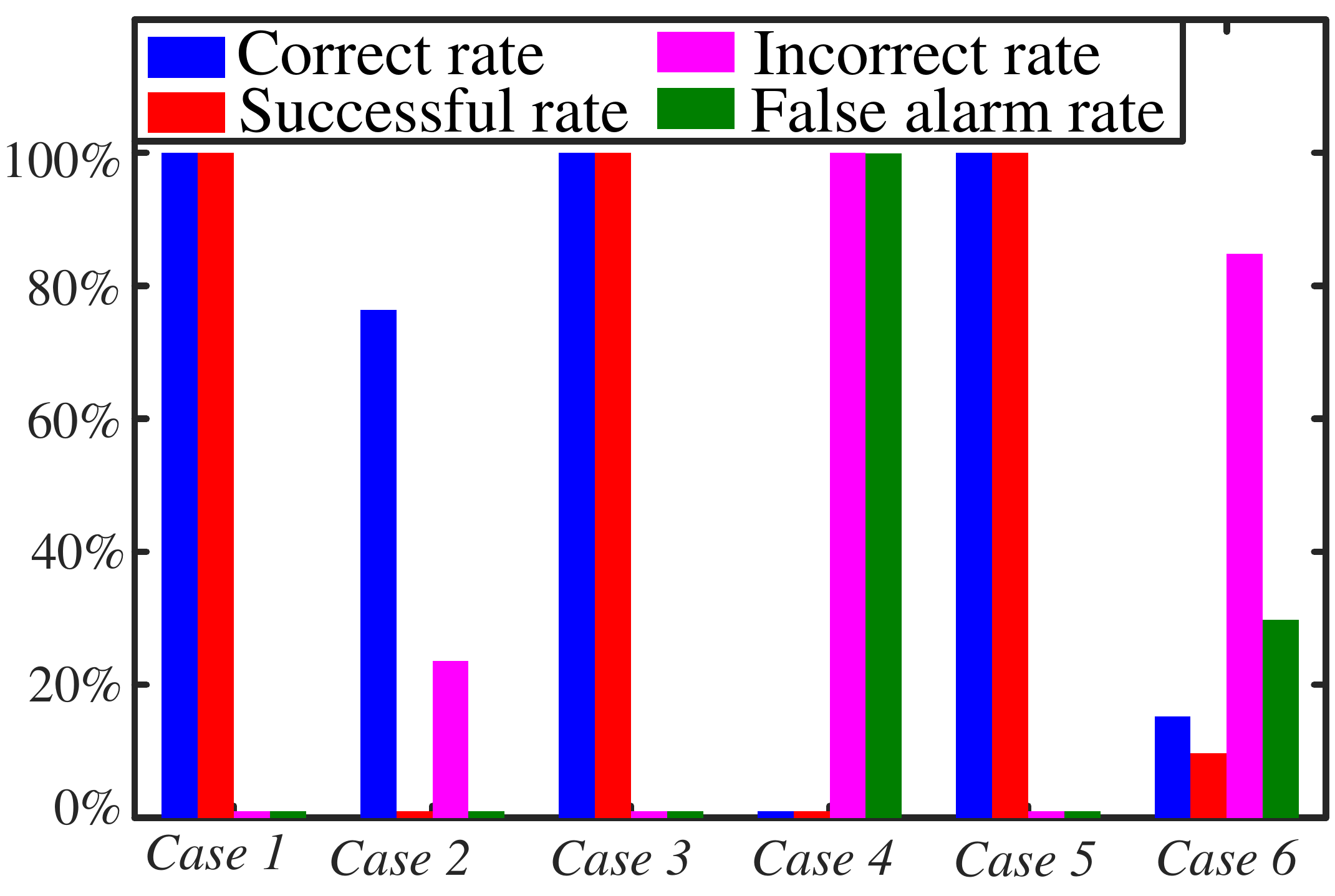} }%
\caption{Correct and successful rate for different cases.}\label{fig:success_flase_alarm}
\end{center}
\end{figure}

If the Lagrange multiplier of the fake outage position has the largest value among others, then this parameter is regarded as a correct attack.
Meanwhile, if the corresponding Lagrange multiplier is larger than the threshold, then the attack is regarded as a successful attack.
By contrast, an incorrect attack is defined as the other Lagrange multiplier or measurement residual of the line, which is not fake outage position, obtains the largest value.
A false alarm is defined when the largest Lagrange multiplier or measurement residual of the line, which is not the fake outage position, is larger than the threshold.
Table \ref{tb:target_lines_selec} shows the selected target lines for the six cases.
Fig. \ref{fig:success_flase_alarm} and Table \ref{tb:target_lines_selec} present the results.

\begin{table} 
\begin{center}
\caption{The parameter setting and the average normalized Lagrange multiplier of correct attack for six cases}\label{tb:target_lines_selec}
\begin{tabular}{|c|c|c|c|c|c|c|c|c|c|c|}
   \cline{1-5}\cline{7-7}
					& $l_{o}$  				& $l_{m}$  		& $\left|\overline{\mathcal{N}}\right|$ & $\left|\overline{\mathcal{E}}\right|$ & & Average $\lambda_{l_{m}}^{N} $ \\ \cline{1-5}\cline{7-7}

   {\it Case $1$}  & ${\bl 13}$    			& ${\rl 17}$  & $4$ 	& $4$		& & $\phantom{0}3.4159$     \\ \cline{1-5}\cline{7-7}
   {\it Case $2$}  & ${\bl 13}$	    		& ${\rl 19}$  & $3$ 	& $2$       & & $\phantom{0}0.1412$     \\ \cline{1-5}\cline{7-7}
   {\it Case $3$}  & ${\bl 21}$         	& ${\rl 23}$  & $10$ 	& $14$      & & $24.9058$     			\\ \cline{1-5}\cline{7-7}
   {\it Case $4$}  & ${\bl 21}$	    		& ${\rl 20}$  & $9$ 	& $11$      & & $0$     				\\ \cline{1-5}\cline{7-7}
   {\it Case $5$}  & ${\bl \phantom{0}4}$   & ${\rl 21}$  & $16$ 	& $19$      & & $10.8885$     			\\ \cline{1-5}\cline{7-7}
   {\it Case $6$}  & ${\bl \phantom{0}4}$	& ${\rl \phantom{0}2}$  & $5$ & $4$ & & $\phantom{0}3.2296$     \\ \cline{1-5}\cline{7-7}
\end{tabular}
\end{center}
\end{table}

Based on the simulation results presented in Fig. \ref{fig:success_flase_alarm}, our proposed selection framework allows a fake outage event to occur every time for the 14-bus, 24-bus, and 118-bus systems.
However, we cannot guarantee the performance of our framework if the target lines are randomly  selected.
The random selection method also leads to a high false alarm rate for the system.
For {\it Case $2$}, the fake outage position can be captured by the operator.
However, the corresponding normalized Lagrange multipliers are all lower than the threshold as shown in Table \ref{tb:target_lines_selec}, thereby leading to a $0 \%$ success rate.
In the 24-bus system, an unexpected fake outage event is revealed to the operator for {\it Case $4$}.
Therefore, if the target lines are randomly selected, then the unpredictable events are revealed to the operator.
In this case, the results cannot be controlled.
For {\it Case $6$}, if the estimation results are not influenced by noise, then line $1$ shows the largest normalized parameter error.
However, this error is below the detection threshold.
Given that this normalized Lagrange multiplier can only exceed the threshold under the influence of noise, and then the noise may also make the errors of other positions exceed the threshold, thereby hiding the expected fake outage information from the operator.
Moreover, given that the observed outage position is not fixed, the operator can easily notice that the system has been injected with false data by the attacker.

\subsection{Comparison of the AC and DC models}

\begin{table} \small
\begin{center}
\caption{The power flow modification for 14-bus system}\label{tb:tran_loss}
\begin{tabular}{|c|c|c|c|c|c|c|c|c|}
\hline
Line number   	& \multicolumn{2}{c|}{Transmission loss (MVA)	} 		  \\
\hline
	    & Calculated 				& Real 	\\
\hline
$12$    & $0.3584 + 0.7460i$      & $0.0736 + 0.1531i$      \\
\hline
${\bl 13}$    & $0.2484 + 0.4892i$      & $0.2484 + 0.4892i$     \\
\hline
$19$    & $0.0175 + 0.0159i$      & $0.0118 + 0.0107i$         \\
\hline
$20$    & $0.2154 + 0.4386i$      & $0.2154 + 0.4386i$        \\
\hline

\end{tabular}
\end{center}
\end{table}

In Seciton \ref{subsec:compare_random_selection}, the attacks designed by AC model can cause very high impact on the system operation.
By contrast, according to \cite{2013-rah-DC-AV-diff}, the attacks constructed by DC model can cause high residual in AC state estimation, such that the attacks can be easily observed.
This is because only voltage phases, real powers, and loads are modifiable in the DC model; however, AC state estimation further needs voltage magnitudes and reactive powers.
Also, the transmission loss has to be considered in the AC state estimation, such that the power flows at the  ``from'' and ``to'' ends are different.
This approach cannot be obtained with DC model.

Given that transmission loss is considered in the formulation, we then use 14-bus system as an example to explain the rationale behind our use of the AC model.
Table \ref{tb:tran_loss} summarizes the transmission losses in the attack region.
The calculated losses are based on $Z_{l} C_{l} $ as shown in (\ref{eq:AC_power_injec}), and the real losses are taken from MATPOWER according to the modified voltage and current system topology.
Based on the results, the real and calculated transmission losses of lines $13$ and $20$ are the same, while those of lines $12$ and $19$ are slightly different.
These results can be attributed to (\ref{eq:Dist_model_flow_convex}) where the second-order cone programming relaxation is applied, and to the fact that the conditions of the boundary buses are bounded in  (\ref{eq:AC_voltage_boundary}).
Therefore, the $\left|\overline{S}_{l} \right|^{2}$ of some lines are not equal to $ W_{j} C_{l}$, and hence making the calculated losses larger than the real losses.
However, these differences are smaller than the power flow.
To measure the influence of miscalculated loss on the dispatch in the attack region, we utilize the following parameter:
\begin{equation}
O_{\rm inf} =  \frac{1}{ \left|\overline{\mathcal{E}}\right|} \sum_{l \in \left|\overline{\mathcal{E}}\right| } \frac{ \Re(Loss_{l}^{\rm cal} - Loss_{l}^{\rm real}) }{ \Re( \overline{\qs}) }
\end{equation}
where $O_{\rm inf}$ is the influence factor, $Loss_{l}^{\rm cal}$ is the calculated  loss, and $Loss_{l}^{\rm real}$ is the real loss of line $l$.
Then, the influence factors for the 14-bus, 24-bus, and 118-bus systems are $0.95\%$, $2.32\%$, and $1.71\%$, respectively.
The comparison results show that the error of the calculation is very small compared to the line flow.
Therefore, the influence is limited, and the calculated error can be regarded as noise when the measurements are entered into the detection mechanism.

In this work, we aim to create an initial failure through the proposed method, and then induce the cascade failure to the power system.
From Section \ref{subsec:14bus_implement} to \ref{subsec:compare_random_selection}, we have already demonstrated that the fake outage line can be continuously appeared as an outage event, and then mislead the control center.
With this approach, an initial failure can be created.
Then, the cascade failure can be propagated from the initial failure as described in \cite{2008-ren-casecade}.
One more thing that we can discuss is the impact after the attack such as economic impact.
Although the impacts due to cascading failure or blackout are the main purpose of the proposed attack strategy, these impacts are hard to be quantized and evaluated.
Therefore, we just discuss the economic impact in terms of the operation cost in Appendix \ref{subsec:compare_cost}.

\section{Conclusion}\label{sec:conclusion}
In this paper, we propose a joint line-removing and line-maintaining attack strategy based on the AC model in which the attacker maliciously injects false data in the cyber layer to cover a physical event in the power system.
The target of the attack strategy is to create an initial failure and then induce the cascade failure in the system.
When launching the attack, the target lines are identified based on the LODF matrix.
The attack region is obtained through the method developed by the BFS algorithm, and then we modify the measurements from the power flow equations.
The simulation results reveal that our proposed scheme successfully misleads the control center and masks the line-outage event.

The potential countermeasures can be separated into two parts.
One possible method is to use historical data to predict the future income.
Then, if the future states have huge difference compared to predicted states, it can be assumed that there is an attack event in the system.
Another method is to study the statistical characteristic of the states in the system.
More specifically, the statistical characteristic of data generated by the attack is different from previous time slots; therefore, the method such as change point detection can be applied.

{\renewcommand{\baselinestretch}{1}
\begin{footnotesize}
\bibliographystyle{IEEEtran}
\bibliography{References_cyber_physical}
\end{footnotesize}}

\appendix

\subsection{Proof of Selection Criteria in Section III-B}\label{subsec:proof_line_selection}
The formulation in (\ref{eq:SE_WLS}) can be rewritten as follows using the Lagrange method :
\begin{equation}
\min_{\widehat{\qv}, P_{e}} ~~~~~   \qr^{T} \qR^{-1} \qr + \qlambda^{T} P_{e}
\end{equation}
The first-order necessary condition of optimality must satisfy the following:
\begin{equation}
\frac{\partial \qr^{T} \qR^{-1} \qr + \qlambda^{T} P_{e}  }{ \partial P_{e} }  = \qH_{p}\qR^{-1}\qr + \qlambda = 0,
\end{equation}
where $\qH_{p}$ is the Jacobian matrix of the measurement functions, $h\left(P_{e}, \mathcal{G},   \widehat{\qv}\right)$, with respect to the network parameter error vector, $P_{e}$.
When SE converges, the Lagrange multiplier vector $\qlambda$ must be expressed as follows:
\begin{equation}
\qlambda = - \qH_{p}\qR^{-1}\qr .
\end{equation}
For the system, $\qR$ is determined, so $\qlambda$ is related to $\qH_{p}$ and $\qr$.
The $\qH_{p}$ is related to the topology and system states, and $\qr$ is caused by the estimation results.
Therefore, the residual vector, $\qr$, has an important influence on the results of bad data and parameter error detection.
When selecting the target lines, if the outage event of one line has a greater influence on the system compared with the outage event of other lines, then the residual can be increased.
For this purpose, we apply the LODF matrix to determine the impact when the specific line is disconnected.

\subsection{Proof of the Efficiency of the Proposed Attack}\label{subsec:proof_effect}
According to \cite{2006-Abur-parameter-SE,2016-Abur-parameter-SE,abur-book}, the normalized $a$-th measurement residual and the corresponding Lagrange multipliers for the $l$-th line can be represented as
\begin{subequations}
\begin{align}
& \left\lvert \lambda_{l}^{N} \right\rvert  = \frac{ \left\lvert \Lambda_{\{ l, l \}} P_{e \{ l \}}  + H_{p \{:, l \}}^{T} R_{ \{ a, a \} }^{-1} S_{ \{ :, a \} }n_{ a  } \right\rvert }{\sqrt{\Lambda_{ \{ l, l \} }}}, \label{eq:nor_para}\\
& \left\lvert r_{a}^{N} \right\rvert  = \frac{ \left\lvert S_{\{ a, a \} } n_{ a }  +  S_{ \{ a, : \} }   H_{p \{:, l \}} P_{e \{ l \}} \right\rvert }{\sqrt{ S_{ \{ a, a \} } R_{ \{ a, a \}}  }}, \label{eq:nor_resi}
\end{align}
\end{subequations}
where $\Lambda_{\{ l, l \}}  $ denotes the $l$-th column and $l$-th row of the parameter covariance matrix, $\mathbf{\Lambda}$, $ S_{\{ a, a\} }$ is the $a$-th column and $a$-th row of the parameter sensitivity matrix, $\qS$, $ R_{\{ a, a\} }$ is the $a$-th column and $a$-th row of the noise covariance matrix, $\qR$, $H_{p \{:, l \}}$ represents the $l$-th column of $\qH_{p}$, $n_{a}$ is the noise interference of $a$-th measurement, and $P_{e \{ l \}}$ is the parameter error of the $l$-th line.

In this section, we explain how the fake outage position obtains the largest value of the normalized Lagrange multiplier compared with the other lines.
The influence of noise is ignored in the derivation.
If the $m$-th line is the fake outage position and has an erroneous parameter, which means  $P_{e \{ m \}} \neq 0, P_{e \{ l \}} =0 (l \neq m)$, then all measurements are correct.
Based on (\ref{eq:nor_para}), we obtain
\begin{subequations}
\begin{align}
& \left| \lambda_{m}^{N} \right|   = \sqrt{ \Lambda_{ \{ m, m \} }}  \left|P_{e, \{m\} } \right|, \\
& \left| \lambda_{l}^{N} \right|   = \frac{ \Lambda_{ \{ m, l \} } }{ \sqrt{ \Lambda_{ \{ l, l \} } } } \left|P_{e, \{m\} } \right|.
\end{align}
\end{subequations}
The ratio of $\left| \lambda_{l}^{N} \right|$ to $\left| \lambda_{m}^{N} \right|$ can be represented as follows:
\begin{equation}\label{eq:para_ratio}
\left| \frac{ \lambda_{l}^{N}  }{  \lambda_{m}^{N} } \right|  = \frac{ \left|\Lambda_{ \{ m, l \} }\right|  }{ \sqrt{  \Lambda_{ \{ l, l \} } \Lambda_{ \{ m, m \} }  } }.
\end{equation}
To show that $\left| \lambda_{m}^{N} \right| $ is larger than the other normalized Lagrange multipliers, we must prove that the ratio in (\ref{eq:para_ratio}) is lower than $1$.
Consider an expectation of a square value that must be positive as
\begin{equation}\label{eq:expect_para}
E  \left\{  \left[   \left( \lambda_{l} - E \left( \lambda_{l} \right) \right)  -  \left( \lambda_{m} - E \left(\lambda_{m} \right) \right) \right]^{2} \right\} \geq 0.
\end{equation}
From the definition of the covariance matrix, we obtain
\begin{subequations}
\begin{align}
& \Lambda_{ \{ m, m \} }  = E \left\{  \left[  \lambda_{m} - E(\lambda_{m}) \right]^{2}  \right\} , \\
& \Lambda_{ \{ l, l \} }  ~~ \,= E \left\{  \left[  \lambda_{l} - E(\lambda_{l}) \right]^{2}  \right\} , \\
& \Lambda_{ \{ m, l \} }  ~\,= E \left\{  \left[  \lambda_{m} - E(\lambda_{m}) \right]  \left[  \lambda_{l} - E(\lambda_{l}) \right]^{\phantom{0}}   \right\}.
\end{align}
\end{subequations}
Eq.\,(\ref{eq:expect_para}) can then be rewritten as follows:
\begin{equation}\label{eq:expect_cov_para}
 \Lambda_{ \{ l, l \} }  - 2\Lambda_{ \{ m, l \} } + \Lambda_{ \{ m, m \} }  \geq 0 .
\end{equation}
To obtain (\ref{eq:expect_para}), the determinant of (\ref{eq:expect_cov_para}) must remain non-positive as follows:
\begin{equation}
\left(  2 \Lambda_{ \{ m, l \} } \right)^{2} - 4 \Lambda_{ \{ l, l \} } \Lambda_{ \{ m, m \} } \leq 0,
\end{equation}
which yields
\begin{equation}
\frac{  \Lambda_{ \{ m, l \} }^{2} }{ \Lambda_{ \{ l, l \} } \Lambda_{ \{ m, m \} } } \leq 1.
\end{equation}
After taking the square root for both sides, we obtain $\left| \lambda_{m}^{N} \right| \geq \left| \lambda_{l}^{N} \right| $ according  to (\ref{eq:para_ratio}).
Therefore, the normalized Lagrange multiplier of the fake outage position obtains the largest value.

\subsection{Economic Impact of the Attacks} \label{subsec:compare_cost}
In this section, we discuss the economic impact of the proposed attack strategy.
However, some impacts, such as the failures in the system or detecting the failures in the system, cannot be easily quantified as a number. 
To cope with this situation, we only discuss the cost that we can calculate, and therefore we only compare the operation cost.  
The operation cost considers the scenario that only the true outage position is disconnected in the system.
Specifically, this is the minimum operation cost that the control center can obtain after observing the fake outage event.
By contrast, this is also the operation cost that the attackers can \emph{at least} cause to the system.

The operation cost increases from $8297.73\$$ to $8331.50\$$, which raises $33.77\$$, in 14-bus system. 
In 118-bus system, the cost increases from $129660.70\$$ to $129726.25\$$, which raises $65.55\$$.
According to the simulation, the operation cost increases after the attacks.
Moreover, if the control center decides to redispatch the flows, the cost can go even higher, and then the failures will happen.

The simulations results reveal that the operation cost slightly increases.
This is because we only consider the situation that a true outage line is disconnected and the control center ignores the fake outage event.
On the other hand, operational costs for the grid to diagnose the outage event and the delay in recovering the system are not included in the comparison; all these damage created to the system is something hard for us to compute. 
If we wish to create more economic impact to the system, there are two directions to further interfere the operation of the power grid.
First, we can attempt to attack several pairs of target lines, which means several real and fake outage positions, in the system.
More specifically, we can have several attackers in the system, and they attempt to jointly attack the power grid.
By contrast, we can propose another attack strategy, which mainly focuses on the impact to the operation cost, so that the operation cost can be higher than the results.
The abovementioned points can be also regarded as our future research topics.

\end{document}